\documentclass[12pt]{article}
\usepackage{epsfig}

\title{A Fast Algorithm for Computing Expected Loan Portfolio Tranche Loss in the Gaussian Factor Model.}

\author{ Pavel Okunev\footnote{This work was supported by the Director, Office of Science, Office  of Advanced Scientific Computing Research, of the U.S. Department of Energy under Contract No. DE-AC03-76SF00098.} \footnote{E-mail: pokunev@math.lbl.gov } \\  Department of Mathematics\\LBNL and UC Berkeley\\ Berkeley, CA 94720}

\date{June 19,2005}

\begin{document}
\maketitle

\begin{abstract}
We propose a fast algorithm for computing the expected tranche loss in the Gaussian  factor model  with arbitrary accuracy using Hermite expansions. No assumptions about homogeneity of the portfolio are made.  It is intended as an alternative to the much slower Fourier transform based methods \cite{MD}.
\end{abstract}

\section{The Gaussian Factor Model}

Let us consider a portfolio of $N$ loans. Let the notional of loan $i$ be equal to the fraction $f_i$ of the notional of the whole portfolio. This means that if loan $i$ defaults  and the entire notional of the loan is lost the portfolio loses fraction $f_i$ or $100f_i\%$  of its value. In practice when a loan $i$ defaults a fraction $r_i$ of its notional will be recovered by the creditors.  Thus the actual loss given default (LGD) of loan $i$ is
\begin{equation}
LGD_i=f_i(1-r_i)
\end{equation}
fraction or
\begin{equation}
LGD_i=100f_i(1-r_i)\%
\end{equation}
of the notional of the entire portfolio.

We now describe the Gaussian m-factor model of portfolio losses from default. The model requires a number of input parameters. For each loan $i$ we are give a probability $p_i$ of its default. Also for each  $i$ and each $k=1,\ldots,m$ we are given a number $w_{i,k}$ such that $\sum_{k=1}^m w_{i,k}^2<1$. The number  $w_{i,k}$  is the loading factor of the loan $i$ with respect to factor $k$. Let $\phi_1, \ldots, \phi_m$ and $\phi^i, i=1,\ldots,N$ be independent standard normal random variables. Let $\Phi(x)$ be the cdf of the standard normal distribution.
In our model loan $i$ defaults if 
\begin{equation}
\sum_{k=1}^m w_{i,k}\phi_k+\sqrt{1-\sum_{k=1}^m w_{i,k}^2}\phi^i<\Phi^{-1}(p_i)
\end{equation}
This indeed happens with probability $p_i$.
  The factors $\phi_1,\ldots,\phi_m$ are usually interpreted  as the state of the global economy, the state of the  regional economy, the state of a particular industry and so on. Thus they are  the factors that affect the default behavior of all or at least a large group of loans in the portfolio. The factors $\phi^1,\ldots,\phi^N$ are interpreted as the idiosyncratic risks of the loans in the portfolio.
 
Let $I_i$ be defined by
\begin{equation}
I_i=I_{\{loan \ i\ defaulted\}}
\end{equation}
  We define the random loss caused by the default of loan $i$ as 
\begin{equation}
L_i=f_i(1-r_i)I_i,
\end{equation}
where $r_i$ is the recovery rate of loan $i$.
The total loss of the portfolio is 
\begin{equation}
L= \sum_i L_i
\end{equation}

An important property of the Gaussian factor model is that the $L_i$'s are not independent of each other. Their mutual  dependence is induced by the dependence of each $L_i$ on the common factors $\phi_1,\ldots,\phi_m$. Historical data supports the conclusion that losses due to defaults on different loans are correlated with each other. Historical data can also be used to calibrate the loadings $w_{i,k}$.the $L_i$'s are not independent of each other. Their mutual  dependence is induced by the dependence of each $L_i$ on the common factors $\phi_1,\ldots,\phi_m$. Historical data supports the conclusion that losses due to defaults on different loans are correlated with each other. Historical data can also be used to calibrate the loadings $w_{i,k}$.

 \section{Analytic Approximation to the Joint Distribution of $\phi_1,\ldots,\phi_m$ and $L$}

When the values of the factors $\phi_1,\ldots,\phi_m$ are fixed, the probability of the default of loan $i$ becomes
\begin{equation}
p^i=\Phi^{-1}\left( \frac{p_i-\sum_kw_{i,k}\phi_k}{\sqrt{1-\sum_kw_{i,k}^2}} \right)
\end{equation}

The random losses $L_i$ become conditionally independent Bernoulli variables with the mean  given by 
\begin{equation}
E_{cond}(L_i)=f_i(1-r_i)p^i
\end{equation}
and the variance given by
\begin{equation}
VAR_{cond}(L_i)=f_i^2 (1-r_i)^2p^i(1-p^i)
\end{equation}   

By the Central Limit Theorem the conditional distribution of the portfolio loss $L$ given the values of the factors $\phi_1,\ldots,\phi_m$  can be approximated by the normal distribution with the mean
\begin{equation}
E_{cond}(L)=\sum_i E_{cond}(L_i)
\end{equation}
and the variance

\begin{equation}
VAR_{cond}(L)=\sum_i VAR_{cond}(L_i)
\end{equation} 

Then the joint distribution of the factors $\phi_1,\ldots,\phi_m$ and the portfolio loss $L$ can be approximated by a distribution with density
\begin{equation}
\rho(\phi_1,\ldots,\phi_m,L)=\rho_{G,E_{cond}(L),VAR_{cond}(L)}(L)\prod_{k=1}^m\rho_{G,0,1}(\phi_k),
\end{equation}
where $\rho_{G,m,v}(x)$ stands for the Gaussian density with mean $m$ and variance $v$.

\section{Expected Loss of a Tranche of Loan Portfolio}

Let $0\leq a<b\leq 1$. We define a tranche loss profile $Tl_{a,b}(x)$ by
\begin{equation}
Tl_{a,b}(x)=\frac{min(b-a,max(x-a,0))}{b-a}
\end{equation}
Number $a$ is called the attachment point of a tranche, while $b$ is called the detachment point  of a  tranche. 
The expected loss of a tranche is then
\begin{equation}
TLoss(a,b)=\int Tl_{a,b}(L)\rho(\phi_1,\ldots,\phi_m,L)d\phi_1\ldots\phi_mL
\end{equation}
This can be rewritten as a double integral
\begin{equation}
\label{theformula}
TLoss(a,b)=\int \int Tl_{a,b}(L)\rho_{G,E_{cond}(L),VAR_{cond}(L)}(L)dL \prod_{k=1}^m\rho_{G,0,1}(\phi_k)d\phi_1\ldots\phi_m
\end{equation}
The inside integral with respect to $L$ is easily done analytically because $L$ has a simple normal distribution  for fixed values of the factors $\phi_1,\ldots,\phi_m$. The outside integral has to be computed numerically. However, since it is an integral of a bounded smooth function with respect to m-dimensional Gaussian density, it is one of the simpler integrals to compute numerically.

\section{Numerical Example}
\begin{figure}[htbp]
\caption{Loan Portfolio Tranche Loss in the Gaussian Single Factor Model}
\epsfig{file=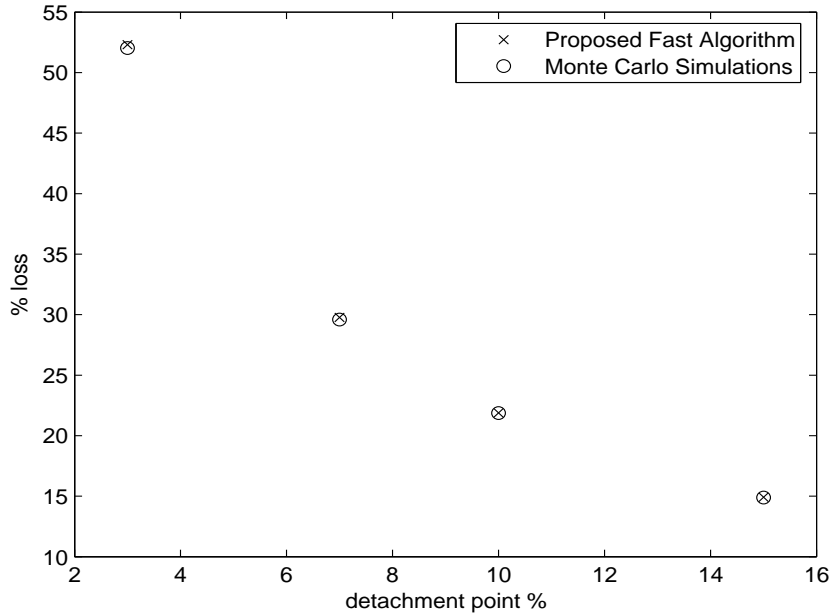,width=5.0 in,height=3.5 in}
\label{resultsfig}
\end{figure}

In this section we apply the proposed algorithm to the single factor Gaussian model of a portfolio with 125 names. We choose a 125 name portfolio because it is the size of the standard DJCDX.NA.IG portfolio.  We choose a single factor model because it is the one most frequently used in practice.  We evaluate the expected loss for four different tranches. All tranches have attachment point $a=0$ or $a=0\%$. The detachment points are 3\%, 7\%, 10\% and 15\%. We take these detachment points because they are the ones most frequently used in practice  in order to evaluate the base correlation. The parameters of the porfolio are
\begin{eqnarray}
f_i&=&\frac{1}{125} \nonumber \\
p_i&=&0.015+\frac{0.05(i-1)}{124} \nonumber \\
r_i&=&0.5-\frac{0.1(i-1)}{124} \nonumber \\
w_{i1}&=&0.5-\frac{0.1(i-1)}{124}
\end{eqnarray}
In Figure \ref{resultsfig} we compare the expected loss computed using $10^6$ Monte Carlo samples with the expected loss computed using formula (\ref{theformula}). The agreement between the two is good.

To obtain the results in Figure \ref{resultsfig} we only needed to perform a single one dimensional numerical integration for each tranche. This is an  improvement over the Moody's FT method \cite{MD} which requires computing a large number of Fourier transforms for each tranche.  Each individual Fourier transform is as computationally expensive as (\ref{theformula}).

\section{The Hermite Expansion of the Conditional Distribution of the Portfolio Loss $L$}

When the approximation to the conditional distribution of $L$ given by the Central Limit Theorem is deemed insufficiently accurate, an arbitrarily  accurate representation of the  conditional distribution of the portfolio loss $L$ can be obtained from its Hermite series expansion. For historical reasons this expansion is also known as the Charlier series expansion \cite{FE}, \cite{CA}.

Let $F(x)$ be the c.d.f. of the  conditional distribution of the portfolio loss $L$. So that 
\begin{equation}
P(L\leq x)=F(x)
\end{equation}
 For  each fixed value of the factors $\phi_1,\ldots,\phi_m$ we define the normalized conditional loss $\tilde{L}$ by
\begin{equation} 
\tilde{L}=\frac{L-E_{cond}(L)}{\sqrt{VAR_{cond}(L)}}
\end{equation}
Let $\tilde{F}(x)$ be the c.d.f. of the  distribution of the normalized   conditional portfolio loss $\tilde{L}$. So that 
\begin{equation}
P(\tilde{L}\leq x)=\tilde{F}(x)
\end{equation}
We define the Hermite polynomial $H_n(x)$ of degree $n$ by
\begin{equation}
H_n(x)=(-1)^ne^{\frac{x^2}{2}}\frac{d^n}{dx^n}e^{\frac{-x^2}{2}}
\end{equation}
Let $c_n$ be defined by
\begin{equation}
c_n=\frac{(-1)^n}{n!} \int_{-\infty}^{\infty}H_n(x)d\tilde{F}(x)
\end{equation}
Then we have 
\begin{equation}
\label{exp}
\tilde{F}(x)=\sum_{i=0}^{\infty} \int_{-\infty}^{x}c_i H_i(t)\frac{e^{\frac{-t^2}{2}}}{\sqrt{2\pi}}dt
\end{equation}
The series above converges in the sense of distributions (generalized functions) \cite{R}. A good reference on the theory of distributions (generalized functions) is \cite{R}.
Let us pick a finite $N$. Then we have
\begin{equation}
\label{expN}
\tilde{F}(x) \approx \sum_{i=0}^{N} c_i \int_{-\infty}^{x}H_i(t)\frac{e^{\frac{-t^2}{2}}}{\sqrt{2\pi}}dt
\end{equation}
As before the approximation is in the sense of  generalized functions. 
Equation (\ref{expN}) implies that the distribution of the normalized   conditional portfolio loss $\tilde{L}$
can be approximated by a distribution with the density
\begin{equation}
\tilde{\rho}(x)=\sum_{i=0}^{N} c_i H_i(x)\frac{e^{\frac{-x^2}{2}}}{\sqrt{2\pi}}
\end{equation}
The function $\tilde{\rho}(x)$ is not necessarily nonnegative and therefore may not be a probability density in the strict sense. However, as is  explained in \cite{R}, this does not affect the validity of our final result (\ref{theformula}). Therefore we may treat $\tilde{\rho}(x)$ as a real probability density.

The distribution of the unnormalized loss $L$ can be approximated by a distribution with density 
\begin{equation}
\label{den1}
\rho(x)=\sum_{i=0}^{N} \frac{c_i}{\sqrt{VAR_{cond}(L)}} H_i\left(\frac{x-E_{cond}(L)}{\sqrt{VAR_{cond}(L)}}\right)\frac{e^{\frac{-\left(\frac{x-E_{cond}(L)}{\sqrt{VAR_{cond}(L)}}\right)^2}{2}}}{\sqrt{2\pi}}
\end{equation}

The joint distribution of the factors $\phi_1,\ldots,\phi_m$ and the portfolio loss $L$ can be approximated by a distribution with density
\begin{equation}
\label{den}
\rho_{joint}(\phi_1,\ldots,\phi_m,L)=\rho(L)\prod_{k=1}^m\rho_{G,0,1}(\phi_k),
\end{equation}
where $\rho_{G,0,1}(x)$ stands for the Gaussian density with mean $0$ and variance $1$.

Observe that the coefficient $c_n$  depends only on the moments of the distribution $\tilde{F}(x)$.
Since $L_i$'s are independent Bernoulli random variables these moments  are known analytically. Thus in the case under consideration all the $c_n$'s are known analytically.

If in equation (\ref{den1}) we set $N=1$ we  obtain  the normal density approximation proposed in the previous section. We show later that it gives good numerical results even when the portfolio size is too small for the normal approximation to be accurate.

\section{Expected Loss of a Tranche of Loan Portfolio}

The expected loss of a tranche  can be written as a double integral
\begin{equation}
\label{theformula2}
TLoss(a,b)=\int \int Tl_{a,b}(L)\rho(L)dL \prod_{k=1}^m\rho_{G,0,1}(\phi_k)d\phi_1\ldots\phi_m
\end{equation}
The inside integral with respect to $L$ can be  done analytically  for fixed values of the factors $\phi_1,\ldots,\phi_m$. The outside integral has to be computed numerically. However, since it is an integral of a bounded smooth function with respect to m-dimensional Gaussian density, it is one of the simpler integrals to compute numerically.

\section{Numerical Example 2}
In this section we test the proposed algorithm on several portfolios of smaller size. For these portfolios the approximation to the conditional distribution of the portfolio loss $L$ given by the Central Limit Theorem is not very accurate, because of their small size.  However, the Hermite expansion produces very good results. We apply the proposed algorithm to the single factor Gaussian model of a portfolio with $n$ names. We take $n$ to be 25 (size of DJ iTraxx Australia), 30  (size of DJ iTraxx ex Japan), 50 (size of DJ iTraxx CJ) and 100 (size of DJCDX.NA.HY).  We choose a single factor model because it is the one most frequently used in practice.  For each $n$ we compute the loss of the equity tranche  with the attachment point $a=0$ or $a=0\%$ and the detachment point 3\%.  The parameters of the porfolio are
\begin{eqnarray}
f_i&=&\frac{1}{n} \nonumber \\
p_i&=&0.015+\frac{0.05(i-1)}{n-1} \nonumber \\
r_i&=&0.5-\frac{0.1(i-1)}{n-1} \nonumber \\
w_{i1}&=&0.5-\frac{0.1(i-1)}{n-1},
\end{eqnarray}
where $i=1,\ldots,n$. Finally, we choose $N=5$ in (\ref{expN}).

\begin{figure}[htbp]
\caption{Equity Tranche Loss in the Gaussian Single Factor Model}
\epsfig{file=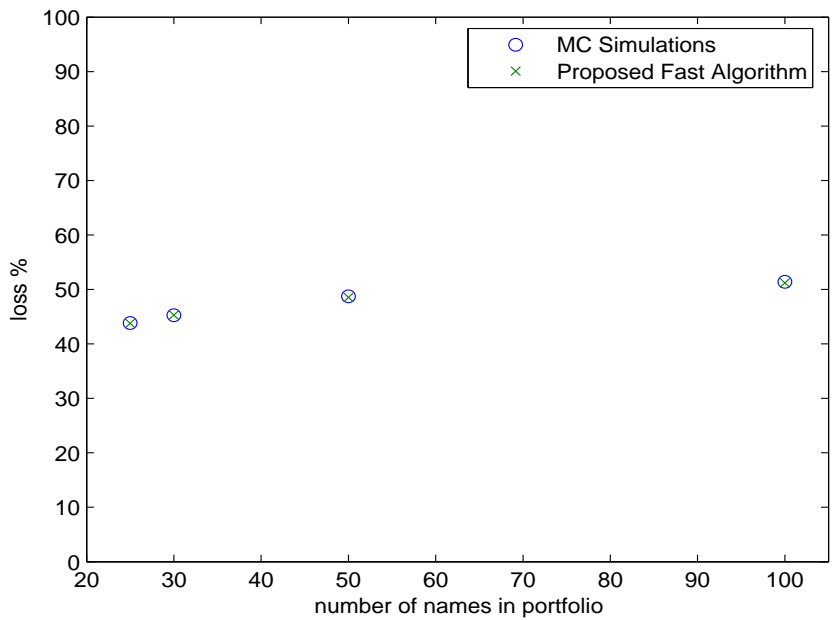,width=5.0 in,height=3.5 in}
\label{resultsfig2}
\end{figure}

In Figure \ref{resultsfig2} we compare the expected loss computed using $10^6$ Monte Carlo samples with the expected loss computed using formula (\ref{theformula2}).\footnote{The author has the code implementing the algorithm described here in MATLAB, VBA for Excel and C.} The agreement between the two is good.

To obtain the results in Figure \ref{resultsfig2} we only needed to perform a single one dimensional numerical integration for each tranche. This is an  improvement over the Fourier transform based methods \cite{MD} which require computing a large number of Fourier transforms for each tranche.  Each individual Fourier transform is as computationally expensive as (\ref{theformula2}).

The Hermite expansion (\ref{expN}) can be used to achieve arbitrary accuracy when the normal approximation is insufficiently accurate. The proposed algorithm is  fast  because the inside integral in (\ref{theformula2}) can  be done analytically.

We also comment that the algorithm can be extended trivially to the case of non-constant recovery rates and recovery rates correlated with the state of the factor variables.

\section{Acknowledgments}
I thank my adviser A. Chorin for his help and guidance during my time in UC Berkeley. I thank Mathilda Regan and Valerie Heatlie for their help in preparing this article. I am also grateful to Ting Lei, Sunita Ganapati and George Wick  for encouraging my interest in finance.  Last, but not least, I thank my family for their constant support.


\begin{thebibliography}{99}

\bibitem{CA} H. Cramer. Mathematical Methods of Statistics.  Princeton University Press, 1954.

\bibitem{MD} A. Debuysscher, M. Szeg\"{o}, M. Freydefront and H. Tabe. Fourier Transform Method-Technical Document. Available from Moody's.

\bibitem{FE} W. Feller. An Introduction to Probability Theory and Its  Applications. Wiley, 1968.

\bibitem{R} M.Reed and B.Simon. Methods of Modern Mathematical Physics, v. 2. Academic press,1972.

\end{thebibliography}
\end{document}